# Quadratic Volume preserving maps

Héctor E. Lomelí
James D. Meiss*
Department of Applied Mathematics
University of Colorado
Boulder, CO 80309

**Abstract**

We study quadratic, volume preserving diffeomorphisms whose inverse is also quadratic. Such maps generalize the Hénon area preserving map and the family of symplectic quadratic maps studied by Moser. In particular, we investigate a family of quadratic volume preserving maps in three space for which we find a normal form and study invariant sets. We also give an alternative proof of a theorem by Moser classifying quadratic symplectic maps.

**AMS classification scheme numbers:**

34C20,34C35,34C37,58F05,70H99

## 1 Introduction

The study of the dynamics of polynomial mappings has a long history both in applied and pure dynamics. For example, such mappings provide simple models of the motion of charged particles through magnetic lenses and are often used in the study of particle accelerators [1]. Moreover, the quadratic, area preserving map, introduced by Hénon [2], is one of the simplest models of chaotic dynamics.

Hénon's study can be generalized in several directions. For example, Moser [3] studied the class of quadratic, symplectic maps, obtaining a useful decomposition and normal form. Here we do the same for more general class of quadratic, volume preserving maps, with one caveat as we discuss below.

Just as symplectic maps arise as Poincaré maps of Hamiltonian flows, volume preserving maps are obtained from incompressible flows, and as such have application to fluid and magnetic field line dynamics [4, 5]. Moreover, one can argue that computational algorithms for flows should obey the "principle of qualitative consistency" [6]: if a flow has some qualitative property then the algorithm should as well. For the case of Hamiltonian flows

*Useful conversations with R. Easton, K. Lenz and B. Peckham are gratefully aknowledge. JDM was supported in part by NSF grant number DMS-9623216.

this leads to the construction of symplectic algorithms. A volume preserving algorithm should be used for a volume preserving flow, such as the motion of passive particle in an incompressible fluid [7, 8, 9, 10].

Some of the properties of symplectic maps generalize to the volume preserving case. For example, a map that is sufficiently close to integrable and nondegenerate in a certain sense has lots of codimension one invariant tori [11, 12, 13], which are absolute barriers to transport [14]. Also, a perturbation of a volume preserving map with a heteroclinic connection can have an exponentially small transversal crossing [15]. Finally, the Birkhoff normal form analysis can be used to study the motion in the neighborhood of fixed points [16]

Another motivation for the study of volume preserving maps is that they can be used as simple models for the study of transport in higher dimensions. The general theory of transport is based on a partition of phase space into regions between which transport is restricted by partial barriers [17]. For example, in two dimensions a partition is formed from intersecting stable and unstable manifolds of a saddle periodic orbit. In higher dimensions an analogous construction requires the existence of codimension one manifolds that separate the space [18]. In most cases it is difficult to find a dynamically natural construction of such manifolds; however, such manifolds do appear in volume preserving maps, and this leads easily to the construction of partial barriers.

The computation and effective visualization of invariant manifolds in higher dimensional maps is itself an interesting problem [19]. In this paper we will study the intersections of the two dimensional stable and unstable manifolds in $\mathbb{R}^3$.

Polynomial maps are also of interest from a mathematical perspective. Much work has been done on the "Cremona maps," that is polynomial maps with constant Jacobians [20]. An interesting mathematical problem concerning such maps is the conjecture proposed by O.T. Keller in 1939:

**Conjecture 1.1 (Real Jacobian Conjecture)** *Let $f : \mathbb{R}^n \to \mathbb{R}^n$ be a Cremona map. Then $f$ is bijective and has a polynomial inverse.*

This conjecture is still open. It is known that injective polynomial maps are automatically surjective and have polynomial inverses [21, 22], so it would suffice to prove that $f$ is injective. It is easy to see (cf. lemma [2.1] below) that for the quadratic case, the condition of volume preservation implies injectivity, thus the Jacobian conjecture holds for quadratic maps.

Even if the conjecture is true, the degree of the inverse of a Cremona map could be large. For example, the upper bound for the degree of the inverse of a quadratic map on $\mathbb{R}^n$ is known to be $2^{n-1}$ [22]. Thus in two dimensions the inverse of a quadratic area preserving mapping is quadratic, as was discussed by Hénon [23, 20]. More generally, Moser showed that quadratic symplectic mappings in any dimension have quadratic inverses [3].

Hénon found the normal form for the quadratic Cremona mapping in the plane. In this paper, we will correspondingly find the normal form for higher dimensional cases, but we assume that the quadratic, volume preserving mapping has a quadratic inverse (it is a "quadratic automorphism"). We give a complete classification of these diffeomorphisms.

Such maps can be written as the composition of an affine volume preserving map and a "quadratic shear." We give necessary and sufficient conditions for such shears to have a quadratic inverse. As a first application of this concept, we give a simple proof of the theorem of Moser [3] for the symplectic case.

We also show that the quadratic automorphism in $\mathbb{R}^3$ can be reduced to one of three normal forms. The generic case has four parameters: two govern the linearization of the map–the trace and second trace of the Jacobian matrix. The remaining two parameters determine the nonlinear terms represented by a single quadratic form in two variables.

## 2   Quadratic Shears

In this section we will study maps of the form

$$x \mapsto x + \frac{1}{2}Q(x)$$

where $Q$ is a vector of quadratic polynomials. Throughout this paper we will write vectors of quadratic polynomials using the form $Q(x) = M(x)x$ where $M : \mathbb{R}^n \to \mathrm{GL}(\mathbb{R}^n)$ is a linear function into the set of $n \times n$ matrices that satisfies the symmetry property $M(x)y = M(y)x$ so that $D_x(M(x)x) = 2M(x)$.

**Definition 2.1** *We say that $f : \mathbb{R}^n \to \mathbb{R}^n$ is a quadratic map in standard form if $f$ can be written as*

$$f(x) = x + \frac{1}{2}M(x)x$$

*where $M$ is a matrix valued linear function that satisfies $M(x)y = M(y)x$.*

It is important to notice that $Df(x) = I + M(x)$.

**Lemma 2.1** *Let $f(x) = x + \frac{1}{2}M(x)x$ be a quadratic map of $\mathbb{R}^n$ in standard form. The following statements are equivalent*

i) *For all $x \in \mathbb{R}^n$, $\det(Df(x)) = 1$.*

ii) *$f$ is bijective with polynomial inverse.*

iii) *$[M(x)]^n = 0$.*

**Proof:**

We will show iii)$\Rightarrow$ii)$\Rightarrow$i)$\Rightarrow$iii).

**iii)$\Rightarrow$ii)**   The condition implies that the matrix $I + M(x)$ is invertible with inverse $I - M(x) + M(x)^2 - \cdots - (-1)^n M(x)^{n-1}$. We notice that we can write

$$f(x) - f(y) = \left(I + M\left(\frac{x+y}{2}\right)\right)(x-y).$$

So the function is injective. Using theorem A in [21], we conclude that $f$ is bijective with a polynomial inverse.

**ii)⇒i)**   In principle, $\det(Df(x))$ and $\det(Df^{-1}(f(x)))$ are polynomials in $x_1, x_2, \ldots, x_n$. However, differentiation of $f^{-1}(f(x)) = x$ gives

$$\det(Df^{-1}(f(x))) \det(Df(x)) = 1,$$

and therefore, since both are polynomials, $\det(Df(x))$ has to be a constant independent of $x$. We notice that $\det(Df(x)) = \det(Df(0)) = \det(I) = 1$.

**i)⇒iii)**   Since $\det(I + M(x)) = 1$ and $M$ is linear in $x$, then for any $\zeta \neq 0$

$$\det(M(x) - \zeta I) = (-1)^n \zeta^n \det(I + M(-\frac{1}{\zeta}x)) = (-1)^n \zeta^n.$$

This implies that the characteristic polynomial of $M(x)$ is $(-\zeta)^n$ and therefore $[M(x)]^n = 0$.
∎

At this point, we restrict to the case of quadratic maps in standard form whose inverse is also quadratic. We will see that the dynamics of such maps is essentially integrable, and is similar to the dynamics of a shear. We first establish the following characterization.

**Lemma 2.2** *Let $f(x) = x + \frac{1}{2}M(x)x$ be a bijective quadratic map of $\mathbb{R}^n$. Then the following statements are equivalent.*

i) $f^{-1}$ *is a quadratic map.*

ii) $M(x)^2 x \equiv 0$.

iii) $M(x)M(y)z + M(y)M(z)x + M(z)M(x)y \equiv 0$.

iv) $f^k(x) = x + \dfrac{k}{2}M(x)x$ *for all $k \in \mathbb{Z}$.*

**Proof:**

We will show i)⇒ii)⇒iii)⇒iv)⇒i).

**i)⇒ii)**   Let

$$f^{-1}(x) = x + \frac{1}{2}N(x)x$$

where $N(x)$ is a matrix valued linear function that satisfies $N(x)y = N(y)x$. Then we have that

$$x = f(f^{-1}(x)) = x + \frac{1}{2}N(x)x + \frac{1}{2}M(x + \frac{1}{2}N(x)x)(x + \frac{1}{2}N(x)x)$$

$$= x + \frac{1}{2}N(x)x + \frac{1}{2}M(x)x + \frac{1}{2}M(x)N(x)x + \frac{1}{8}M(N(x)x)N(x)x$$

This implies that $N(x)x = -M(x)x$ and $M(x)M(x)x = 0$.

**ii)⇒iii)**   By linearity and symmetry of $M$.

**iii)⇒iv)**   Let $g_k(x) = x + \frac{k}{2}M(x)x$. Then

$$g_k(g_l(x)) = x + \frac{l}{2}M(x)x + \frac{k}{2}M(x + \frac{l}{2}M(x)x)(x + \frac{l}{2}M(x)x)$$

$$= x + \frac{l+k}{2}M(x)x + \frac{kl}{2}M(x)M(x)x + \frac{kl^2}{8}M(M(x))M(x)x$$

$$= x + \frac{l+k}{2}M(x)x - \frac{kl^2}{8}[M(x)M(x)M(x)x + M(x)M(M(x))x]$$

$$= x + \frac{l+k}{2}M(x)x.$$

Therefore $g_k \circ g_l = g_{k+l}$. On the other hand $g_1 = f$ and $g_0 = id$. This implies that $g_k = f^k$.

**iv)⇒i)**   Clear. ■

**Definition 2.2** *Let $f : \mathbb{R}^n \to \mathbb{R}^n$ be given, in standard form, by $f(x) = x + \frac{1}{2}M(x)x$. If $f$ satisfies any of the conditions of lemma [2.2], we will say that $f$ is a quadratic shear.*

A simple family of quadratic shears is determined by any vector $v \in \mathbb{R}^n$ and a symmetric matrix $P$ such that $Pv = 0$ according to $M(x)y = (x^T Py)v$, for then

$$M(x)^2 x = (x^T Pv)(x^T Px)v = 0.$$

We will see that, at least in the case $n = 3$, this is the most general quadratic shear. Moser's normal form for symplectic, quadratic maps [3] shows that the higher dimensional case is not quite this simple. From now on, we will concentrate on the special case $n = 3$.

**Theorem 2.1** *A function $f : \mathbb{R}^3 \to \mathbb{R}^3$ is a quadratic shear in $\mathbb{R}^3$ if and only if there is a vector $v \in \mathbb{R}^3$ and a $3 \times 3$ symmetric matrix $P$ such that $Pv = 0$ and*

$$f(x) = x + \frac{1}{2}(x^T Px)v$$

**Proof:**

Since $f$ is a bijection, we can define a new function $g : S^2 \to S^2$. on the unit two dimensional sphere $S^2 \subset \mathbb{R}^3$, in the following way.

$$g(x) = \frac{f(x)}{|f(x)|}.$$

Using standard theorems of algebraic topology [24], we can argue that $g$ has either a fixed point or an antipodal point (a point such that $g(x) = -x$). In any case, there is a constant $K \in \mathbb{R} \setminus \{0\}$ and a vector $x_0 \neq 0$ such that $f(x_0) = Kx_0$. We will show that $K = 1$. We notice that $x_0$ satisfies the following

$$f(x_0) = Kx_0 = x_0 + \frac{1}{2}M(x_0)x_0,$$

$$f^{-1}(Kx_0) = x_0 = Kx_0 - \frac{K^2}{2}M(x_0)x_0,$$

and therefore

$$K^2(Kx_0 - x_0) = \frac{K^2}{2}M(x_0)x_0 = Kx_0 - x_0.$$

We conclude that $K$ satisfies $K^3 - K^2 - K + 1 = 0$ and hence $K = 1$ or $K = -1$. Clearly, since $f$ is a bijection, $f(\frac{1}{2}x_0) \neq 0$ implies that $M(x_0)x_0 \neq -4x_0$. We conclude that $f(x_0) \neq -x_0$ and therefore $K = 1$.

It is clear that $M(x_0)x_0 = 0$. Without loss of generality, we can assume that $x_0 = e_1 = (1, 0, 0)$. Notice that $M(e_1)$ has to have the form

$$M(e_1) = \begin{pmatrix} 0 & \gamma_1 & \gamma_2 \end{pmatrix}$$

where $\gamma_1, \gamma_2 \in \mathbb{R}^3$. This fact, together with lemma [2.2], implies that the matrix $M(e_1)^2 = 0$, and therefore $\gamma_1$ and $\gamma_2$ are parallel. We can perform a linear change of coordinates and assume without loss of generality that $f$ has the form

$$f\begin{pmatrix} x_1 \\ x_2 \\ x_3 \end{pmatrix} = \begin{pmatrix} x_1 \\ x_2 \\ x_3 \end{pmatrix} + x_1 x_2 \begin{pmatrix} \alpha \\ 0 \\ \beta \end{pmatrix} +$$

$$\frac{x_2^2}{2}\begin{pmatrix} \mu_1 \\ \mu_2 \\ \mu_3 \end{pmatrix} + x_2 x_3 \begin{pmatrix} \nu_1 \\ \nu_2 \\ \nu_3 \end{pmatrix} + \frac{x_3^2}{2}\begin{pmatrix} \eta_1 \\ \eta_2 \\ \eta_3 \end{pmatrix}$$

Let $M_1 = M(e_1), M_2 = M(e_2)$ and $M_3 = M(e_3)$ where $e_1 = (1, 0, 0), e_2 = (0, 1, 0)$ and $e_3 = (0, 0, 1)$. It is easy to see that

$$M_1 = \begin{pmatrix} 0 & \alpha & 0 \\ 0 & 0 & 0 \\ 0 & \beta & 0 \end{pmatrix}$$

$$M_2 = \begin{pmatrix} \alpha & \mu_1 & \nu_1 \\ 0 & \mu_2 & \nu_2 \\ \beta & \mu_3 & \nu_3 \end{pmatrix}$$

and

$$M_3 = \begin{pmatrix} 0 & \nu_1 & \eta_1 \\ 0 & \nu_2 & \eta_2 \\ 0 & \nu_3 & \eta_3 \end{pmatrix}.$$

To finish the proof, we need to show that the column vectors $\mu, \nu, \eta$ and $(\alpha, 0, \beta)$ of $M_1, M_2, M_3$ are parallel to each other. We will show step by step that several of the entries are zero. We have two cases.

- $\beta \neq 0$. Using lemma [2.2] we conclude that $2M_3^2 e_1 + M_1 M_3 e_3 = 0$. This implies that $\eta_2 = 0$. We also have that $M_3^3 = 0$, so $\nu_2 = 0$ and $\eta_3 = 0$. The condition $M_2^2 e_2 = 0$ implies that $\mu_2 = 0$ and this together with $M_2^3 = 0$ implies that $\nu_3 = -\alpha$. Using the equation $M_2 M_3 e_3 + 2M_3^2 e_2 = 0$ we find that $\eta_1 = 0$. Using $M_2^3 = 0$ and $M_2^2 e_2 = 0$, we find that the column vectors of $M_2$ are parallel, and the rest is clear.

- $\beta = 0$. The condition $M_2^3 = 0$ implies that $\alpha = 0$, $\nu_3 = -\mu_2$ and $\mu_2^2 + \nu_2\mu_3 = 0$. The condition $M_3^3 = 0$ implies that $\eta_3 = -\nu_2$ and $\nu_2^2 - \mu_2\eta_2 = 0$. On the other hand, $M_2^2 e_2 = 0$ implies that $\mu_1\mu_2 + \mu_3\nu_1 = 0$ and $M_3^2 e_3 = 0$ implies that $\nu_1\eta_2 - \eta_1\nu_2 = 0$.

  So it is enough to show that $\mu_1\nu_2 - \nu_1\mu_2 = 0$ and $\nu_1\nu_2 - \eta_1\mu_2 = 0$. Clearly, if $\eta_2 = 0$ then $\mu_2 = 0$ and we would be done. So, we can assume that $\eta_2 \neq 0$.

  If $\eta_2 \neq 0$ then $\eta_1\mu_2 = \eta_1\nu_2^2/\eta_2 = \nu_1\nu_2$. If $\nu_2 = 0$ then $\mu_2 = 0$ and we would be done. Assume that $\nu_2 \neq 0$ and $\eta_2 \neq 0$. This implies that $\mu_1\nu_2 = \mu_1\mu_2\eta_2/\nu_2 = -\mu_3\nu_1\eta_2/\nu_2 = \mu_2^2\nu_1\eta_2/\nu_2^2 = \mu_2\nu_1$. ∎

## 3 Quadratic Symplectic Maps

In this section we use the characterization of quadratic shears in lemma [2.2] to give an alternate proof of the result of Moser [3] for quadratic symplectic maps. Recall that a map $f$ is symplectic if $\omega(Dfv, Dfv') = \omega(v, v')$ for all vectors $v, v' \in \mathbb{R}^{2n}$ where $\omega$ is the standard symplectic form $\omega(v, v') = v^T J v'$ and $J$ is the $2n \times 2n$ matrix,

$$J = \begin{pmatrix} 0 & I \\ -I & 0 \end{pmatrix} .$$

**Theorem 3.1** *Let $f$ be a quadratic symplectic map of $(\mathbb{R}^{2n}, \omega)$. Then $f$ can be decomposed as $f = T \circ S$ where $T$ is affine symplectic and $S$ is a symplectic quadratic shear. Furthermore, if $S$ is any symplectic quadratic shear, then there is a symplectic linear map $\lambda$ such that $\lambda \circ S \circ \lambda^{-1}(q, p) = (q + \nabla V(p), p)$.*

**Proof:**

Let $b = f(0)$ and $L = Df(0)$. Clearly $L$ is a symplectic matrix and if we let $T(x) = Lx + b$ then $S = T^{-1} \circ f$ is a symplectic quadratic map in standard form, i.e. $S(x) = x + \frac{1}{2}M(x)x$, where $M(x)$ is linear in $x$ and satisfies the symmetry property $M(x)y = M(y)x$. Then $S$ is symplectic providing

$$(I + M(x))^T J (I + M(x)) = J .$$

Homogeneity of $M(x)$ implies that

$$M(x)^T J = J^T M(x) , \tag{1}$$

and

$$M(x)^T J M(x) = 0 . \tag{2}$$

Using (1) in (2) gives $0 = M(x)^T J M(x) = J^T M(x) M(x)$, and since $J$ is nonsingular this implies

$$M(x)^2 = 0 . \tag{3}$$

Then lemma [2.2] implies that $S$ is a quadratic shear.

To finish the proof, we follow Moser [3] and define the null space of $M(x)$ in the following way
$$N = N(M) = \{y \in \mathbb{R}^{2n} : M(x)y = 0, \forall x \in \mathbb{R}^{2n}\} = \{y \in \mathbb{R}^{2n} : M(y) = 0\} \ .$$
Recall [25] that the $\omega$−orthogonal complement of a subspace $E \subset \mathbb{R}^{2n}$ is defined by $E^\perp = \{v \in \mathbb{R}^{2n} : \omega(v, v') = 0, \forall v' \in E\}$. We will show that $N^\perp \subset N$. For that purpose, we will use the following fact
$$M(z)M(x)y = M(x-y)^2 z = 0 \tag{4}$$
that follows from lemma [2.2], linearity and symmetry.

Let $u \in N^\perp$ and $x \in \mathbb{R}^{2n}$. Now for any $y \in \mathbb{R}^{2n}$, (4) implies that $M(x)y \in N$. Therefore $\omega(y, M(x)u) = y^T J M(x) u = -y^T M(x)^T J u = -\omega(M(x)y, u) = 0$. This implies that $M(x)u = 0$ and hence $u \in N$. Standard theorems in symplectic geometry (cf. [25]) imply that it is possible to find a Lagrangian space $F$ such that $N^\perp \subset F^\perp = F \subset N$ and a symplectic linear transformation $\lambda$ such that
$$\lambda(F) = \{(q, p) \in \mathbb{R}^n \times \mathbb{R}^n : p = 0\}.$$
Clearly, if $S(x) = I + \frac{1}{2}M(x)x$ is a symplectic quadratic shear, then so is $\tilde{S} = \lambda \circ S \circ \lambda^{-1}$. Assume that $\tilde{S}(x) = I + \frac{1}{2}\tilde{M}(x)x$. Then $\lambda(F) \subset N(\tilde{M})$. Therefore for all $(q, p) \in \mathbb{R}^n \times \mathbb{R}^n$,
$$\tilde{M}(q,p)(q,p) = \tilde{M}(q,p)(0,p) = \tilde{M}(0,p)(q,p) = \tilde{M}(0,p)(0,p) \ .$$
Since, in general, the matrix $\tilde{M}(0, p)$ can be written in $n \times n$ blocks as
$$\tilde{M}(0, p) = \begin{pmatrix} A(p) & B(p) \\ C(p) & D(p) \end{pmatrix},$$
then $\tilde{M}(0,p)(q,0) = 0$ implies $A(p) = C(p) = 0$. Moreover, (1) implies $D(p) = 0$ and $B(p)^T = B(p)$. Thus finally we see that
$$\tilde{M}(q,p)(q,p) = (B(p)p, 0)$$
where $B(p)p$ is a gradient vector field. ∎

# 4 Normal Form in $\mathbb{R}^3$

In this section we give normal forms for a quadratic diffeomorphism $f$ of $\mathbb{R}^3$ that preserves volume and has a quadratic inverse. Now lemma [2.2] implies that if we let $b = f(0)$ and $L = Df(0)$, and $T(x) = Lx + b$, then the map $S = T^{-1} \circ f$ is a quadratic shear. Then theorem [2.1] implies that $S$ is of the form $S(x) = x + \frac{1}{2}(x^T P x)v$ where $v \in \mathbb{R}^3$ and $P$ is a symmetric matrix such that $Pv = 0$. Depending on the relation between $L$ and $v$, we have three cases possible cases; these can by distinguished by considering the space
$$Z(v, L) = span\{v, Lv, L^2 v\} \ .$$

**Theorem 4.1** *Let $f : \mathbb{R}^3 \to \mathbb{R}^3$ be a quadratic volume preserving diffeomorphism. Then $f$ can be written as the composition of an affine map $T$ and a quadratic shear $S$, $f = T \circ S$, where $S(x) = x + \frac{1}{2}(x^T P x)v$, $v \in \mathbb{R}^3$ and $P$ is a symmetric matrix such that $Pv = 0$. Moreover, $f$ is affinely conjugate to one of three possible normal forms, depending on the dimension of $Z(v, L)$:*

i) $\dim Z(v, L) = 3$. *The map $f$ is conjugate to*

$$\begin{pmatrix} \alpha + \tau x - \sigma y + z + Q(x, y) \\ x \\ y \end{pmatrix} \qquad (5)$$

*where $\tau$ and $\sigma$ are the trace and second trace of $L$, and $Q(x, y) = ax^2 + bxy + cy^2$ is a quadratic form.*

ii) $\dim Z(v, L) = 2$. *The map $f$ is conjugate to*

$$\begin{pmatrix} x_0 + \alpha x + y + Q(x, z) \\ y_0 - \beta x \\ z_0 + \frac{1}{\beta} z \end{pmatrix}.$$

iii) $\dim Z(v, L) = 1$. *The map $f$ is conjugate to*

$$\begin{pmatrix} x_0 + \alpha x + Q(y, z) \\ y_0 - \frac{1}{\alpha} z \\ z_0 + y + \beta z \end{pmatrix}.$$

**Proof:**

We know that $f = L(x + \frac{1}{2}(x^T P x)x) + b$, and $Pv = 0$. To obtain the first normal form, perform a linear change of coordinates, $x = U\xi$. Since the vectors $v, Lv$, and $L^2v$ are linearly independent, the transformation $U$ can be defined by the following equations

$$\begin{aligned} U^{-1}v &= e_3 & Ue_3 &= v \\ U^{-1}Lv &= e_1 & Ue_1 &= Lv \\ U^{-1}L^2v &= e_2 + \tau e_1 & Ue_2 &= L^2v - \tau Lv \end{aligned}$$

where, as we will see below, we will choose $\tau = Tr(L)$. In the new coordinates the map becomes

$$\begin{aligned} \xi' &= U^{-1} f(U(\xi)) \\ &= U^{-1}b + U^{-1}LU\xi + \frac{1}{2}\left(\xi^T U^T P U \xi\right) U^{-1} Lv \\ &= \xi_o + U^{-1}LU\xi + e_1 \tilde{Q}(\xi, \xi) \end{aligned}$$

where $\tilde{Q}(\xi_1, \xi_2) = \frac{1}{2}\left(\xi_1^T U^T P U \xi_2\right)$. Note that $\tilde{Q}(\xi, e_3) = \frac{1}{2}\left(\xi^T U^T Pv\right) = 0$, so in the new coordinates the quadratic terms depend only on the first and second components. Moreover in this coordinate system

$$U^{-1}LUe_1 = U^{-1}L^2v = e_2 + \tau e_1$$
$$U^{-1}LUe_2 = U^{-1}\left(L^3v - \tau L^2v\right)$$
$$U^{-1}LUe_3 = U^{-1}Lv = e_1$$

The second equation can be simplified by noting that the characteristic equation for the matrix $L$ is satisfied by $L$ itself, and so $L^3 - \tau L^2 + \sigma L - I = 0$, where $\tau = Tr(L)$ and $\sigma = Tr_2(L)$, the "second trace" of the matrix $L$, thus we get $U^{-1}LUe_2 = U^{-1}(I - \sigma L)v = e_3 - \sigma e_1$. Thus we obtain

$$U^{-1}LU = \begin{pmatrix} \tau & -\sigma & 1 \\ 1 & 0 & 0 \\ 0 & 1 & 0 \end{pmatrix}$$

Upon reverting to $(x, y, z)$ as the names for the coordinates we get

$$U^{-1}f(U(x)) = \begin{pmatrix} x_0 \\ y_0 \\ z_0 \end{pmatrix} + \begin{pmatrix} \tau x - \sigma y + z + Q(x,y) \\ x \\ y \end{pmatrix}$$

To simplify this map further, we can conjugate, using the translation

$$(x, y, z) \mapsto (x, y + y_0, z + y_0 + z_0),$$

to a map with $x_0 = \alpha$, $y_0 = 0$ and $z_0 = 0$. This is the promised form.

For the second case, assume that $L^2v = \alpha Lv - \beta v$, for some nonzero $\alpha$ and $\beta$. This implies that the characteristic polynomial for $L$ factors as $(L - 1/\beta I)(L^2 - \alpha L + \beta I) = 0$, and therefore, since $L$ is nondegenerate, there exists a vector $w \notin Z(v, L)$ such that $Lw = \frac{1}{\beta}w$. We define the following change of coordinates.

$$\begin{aligned} U^{-1}v &= e_2 & Ue_2 &= v \\ U^{-1}Lv &= e_1 & Ue_1 &= Lv \\ U^{-1}w &= e_3 & Ue_3 &= w \end{aligned}$$

As before, we note that in the new coordinates the quadratic term satisfies $\tilde{Q}(e_2, \xi) = 0$, so in the new coordinates the quadratic terms depend only on the first and third components. Moreover in this coordinate system we obtain

$$U^{-1}LU = \begin{pmatrix} \alpha & 1 & 0 \\ -\beta & 0 & 0 \\ 0 & 0 & \frac{1}{\beta} \end{pmatrix}$$

This implies the form for the second case.

For the third case, assume that $Lv = \alpha v$. We notice that there exist a vector $w \notin Z(v, L)$ such that $Z(w, L) \oplus Z(v, L) = \mathbb{R}^3$. In fact, we can also find a constant $\beta$ such that $L^2 w - \beta L w + \frac{1}{\alpha} w = 0$. We define the following change of coordinates.

$$U^{-1} v = e_1 \qquad\qquad U e_1 = v$$
$$U^{-1} w = e_2 \qquad\qquad U e_2 = w$$
$$U^{-1} L w = e_3 \qquad\qquad U e_3 = L w$$

As before, we note that in the new coordinates the quadratic term is $\tilde{Q}(e_1, \xi) = 0$, so in the new coordinates the quadratic terms depend only on the second and third components. Moreover in this coordinate system we obtain

$$U^{-1} L U = \begin{pmatrix} \alpha & 0 & 0 \\ 0 & 0 & -\frac{1}{\alpha} \\ 0 & 1 & \beta \end{pmatrix}$$

This implies the form for the last case. ∎

## 5  Dynamics

The dynamics of the second and third cases of theorem [4.1] are essentially trivial. In case ii), the $z$ dynamics decouples from the $(x,y)$ dynamics. There are four special cases:

- $|\beta| \neq 1$. The plane $z = \beta z_0 / (\beta - 1)$ is invariant, and is either a global attractor ($|\beta| > 1$) or repellor ($|\beta| < 1$). On the plane the dynamics is linear.

- $\beta = 1$, $z_0 \neq 0$. All orbits are unbounded.

- $\beta = 1$, $z_0 = 0$. Every plane $z = c$ is invariant, and the dynamics reduces to a two dimensional area preserving Hénon map on each plane.

- $\beta = -1$. Each plane $z = c$ is fixed under $f^2$. Restricted to a plane, $f^2$ is the composition of two orientation reversing Hénon maps.

For case iii) the $(y, z)$ dynamics is linear and decouples from the $x$ dynamics. Generically, there is an invariant line on which the dynamics is affine. The invariant line can have any stability type.

### 5.1  Generic case

Equation (5) is the only nontrivial case. In general this map has six parameters, one from the shift, two from the linear matrix (the two coefficients of its characteristic polynomial) and the three coefficients of $Q$. However, generically, two of these parameters can be eliminated.

Write the quadratic form as $Q(x, y) = ax^2 + bxy + cy^2$. Generically $a + b + c \neq 0$ and we can we can apply a scaling transformation to set $a + b + c = 1$. Similarly if $b + 2c \neq 0$ the parameter $\sigma$ can be eliminated using the diagonal translation

$$(x, y, z) \mapsto (x + \gamma, y + \gamma, z + \gamma), \quad \gamma = \sigma / (b + 2c)$$

In this way, we obtain the final, generic form

$$\begin{pmatrix} x' \\ y' \\ z' \end{pmatrix} = \begin{pmatrix} \alpha + \tau x + z + ax^2 + bxy + cy^2 \\ x \\ y \end{pmatrix} \qquad a + b + c = 1 \qquad (6)$$

There are four parameters in the system. Even if $a + b + c = 0$ and/or $b + 2c = 0$, then other normalizations can be chosen to eliminate two of the parameters in (5). We will not study these special cases.

## 5.2 Periodic orbits

Generically we can assume that $a + b + c = 1$ for the quadratic form in (5). The map (5) has at most two fixed points

$$x = y = z = x_\pm = \frac{1}{2}\left(-\tau + \sigma \pm \sqrt{(\tau - \sigma)^2 - 4\alpha}\right) \qquad (7)$$

born in a saddle node bifurcation at $(\tau - \sigma)^2 - 4\alpha = 0$. The characteristic polynomial of the linearized map at the fixed points is

$$\lambda^3 - t\lambda^2 + s\lambda - 1 = 0$$

where the trace $t$ and second trace $s$ are

$$\begin{aligned} t_\pm &= \tau + (2a + b)x_\pm \\ s_\pm &= \sigma - (2c + b)x_\pm \end{aligned}$$

We notice that

$$\begin{aligned} t_\pm - s_\pm &= \pm\sqrt{(\tau - \sigma)^2 - 4\alpha} \\ t_\pm - s_\mp &= \pm(a - c)\sqrt{(\tau - \sigma)^2 - 4\alpha}. \end{aligned}$$

The corresponding eigenvalues are illustrated in Fig 1. It is easy to see (using the symmetric polynomials) that there are two lines in $(t, s)$ space where the stability changes: the saddle-node line $t = s$ corresponds to an eigenvalue 1, and the period doubling line $t + s = -2$ corresponds to an eigenvalue $-1$. At the point $t = s = -1$ where they cross the eigenvalues are necessarily $(-1, -1, 1)$. Note also that when $-1 \le t = s \le 3$ there is a pair of eigenvalues on the unit circle. There are two other curves of interest in the stability plane–these correspond to a double eigenvalue $\lambda_1 = \lambda_2 = r$, or

$$2r + 1/r^2 = t \qquad r^2 + 2/r = s$$

This gives the two curves shown in Fig 1. One has a cusp at $t = s = 3$, where we have the triple root $\lambda = 1$. The second crosses the saddle-node and period doubling lines at $t = s = -1$. These are the two codimension two points.

A fixed point with a one dimensional unstable manifold is called *type A* and one with a one dimensional stable manifold is called *type B*. The saddle-node and period doubling lines divide the plane into quadrants which alternate between type $A$ and $B$.

Having a pair of fixed points one of *type A* and one of *type B*, has interesting consequences for our map. For instance, the two dimensional manifolds serve as partial barriers to transport. Generically, they intersect along a one dimensional manifold. We have computed numerically some pairs of two dimensional stable and unstable manifolds. As an example, see figure 2.

We have noticed that varying the parameter makes the one dimensional intersection bifurcate. Further investigation in this direction is the subject of future papers and more complete treatment will appear elsewhere.

The two fixed points (7) are born on the line $t = s$ and move to opposite sides of this line for $(\tau - \sigma)^2 > 4\alpha$ ($x_+$ is always on the right side). If $|a - c|$ is small, they are on the same side of the period doubling line, so that one is type $A$ and the other type $B$; however, when this parameter is large enough they can be on opposite sides, and therefore of the same *type*. This is determined by the sign of

$$s_\pm + t_\pm + 2 = 2 + \tau + \sigma + 2(a - c)x_\pm$$

When $a = c$, we have $t_\pm = s_\mp$ so that the eigenvalues of the two fixed points are reciprocal (see section 5.3 for the explanation of this).

Remember that, generically and without loss of generality, we can assume that $\sigma = 0$. Therefore, we can plot stability diagrams for different values of $\tau$ and $\alpha$. The stability diagram in the $(\tau, \alpha)$ plane for the $a = c$ case is shown in Fig 3. A more general case is shown in Fig 4.

Periodic orbits can be studied by converting the map into a third order difference equation. Let $(x_t, y_t, z_t), t = 0, 1, \ldots$ be a trajectory of the map (5), then the map can be written as

$$x_{t+1} = \alpha + \tau x_t - \sigma x_{t-1} + x_{t-2} + Q(x_t, x_{t-1}) .$$

Now it is clear that if this system of $n$ quadratic equations has a finite number of solutions, then there are at most $2^n$. We can rule out this degeneracy is most cases

**Lemma 5.1** *Suppose a and c are not both zero, and let $\mu_\pm$ be the two (possibly complex) solutions of $Q(\mu, 1) = 0$. Then if $\mu_+^k \mu_-^{n-k} \neq 1$ for some integer $0 \leq k \leq n$, the number of fixed points of $f^n$ for the map (5) is at most $2^n$.*

**Proof:**

The "nonlinear alternative" [26] asserts that the number of complex solutions, counted with multiplicity, of a system of $n$ polynomial equations in $n$ variables is precisely the product of the degrees of the polynomials providing the system of equations obtained by discarding all terms but those of the highest degree in each equation has only the trivial solution. For our case, the resulting system is

$$Q(x_t, x_{t-1}) = ax_t^2 + bx_t x_{t-1} + cx_{t-1}^2 = 0, \ t =, 1 \ldots n, \ x_n = x_0 .$$

If any one of the $x_t = 0$, then they all are zero, unless $a = c = 0$. Otherwise the general nonzero solution to this system is $x_t = \mu_\pm x_{t-1}$, where $a\mu^2 + b\mu + c = 0$. Setting $x_n = x_0$ requires $\mu_+^k \mu_-^l = 1$, where $k + l = n$. ∎

For example, typically there are at most four fixed points of $f^2$, giving a single period two orbit in addition to the two fixed points of $f$. However, there could be infinitely many period two orbits or none when $\mu_+\mu_- = 1$, giving $a = c$ or $\mu_\pm^2 = 1$, giving $b = \pm(a+c)$. As an example, when $a = c = b/2$ and $\sigma + \tau + 2 = 0$, the line $(x, \delta - x, x) \mapsto (\delta - x, x, \delta - x)$ has period two where $\delta$ is defined by

$$\alpha + (1+\sigma)\delta + a\delta^2 = 0$$

## 5.3 Reversibility

A map is reversible if it is conjugate to its inverse by a diffeomorphism $h$ that is an involution, thus

$$h \circ f = f^{-1} \circ h \ , \ h^2 = I \ .$$

Some of the quadratic maps that we have considered have a reversor. Assume the generic case $a + b + c \neq 0$ (or equivalently $Q(1,1) \neq 0$, where $Q$ is the quadratic form given in (5)). It is easy to see that, if $a = c$, then the map (5) has a reversor given by

$$h(x,y,z) = -\begin{pmatrix} z + \eta \\ y + \eta \\ x + \eta \end{pmatrix}.$$

where $\eta = \tau - \sigma/(a+b+c)$.

Note that when $f$ is reversible and has fixed points, then the two fixed points have reciprocal eigenvalues–so if one is type A, the other is type B. Moreover, if the eigenvalues are complex, then the rotation rates have the same magnitudes at the two fixed points.

**Lemma 5.2** *Let $f$ be a quadratic map in normal form (5). Assume, generically, that the quadratic form satisfies $Q(1,1) \neq 0$ and $(\tau - \sigma)^2 \neq 4\alpha Q(1,1)$ Then $f$ is smoothly reversible if and only if $Q(x,y) = Q(y,x)$.*

**Proof:**

Without loss of generality, we assume that $Q(1,1) = a+b+c = 1$ and $(\tau-\sigma)^2 - 4\alpha \neq 0$. Extend the map to $\mathbb{C}^3$. The imposed conditions imply that the map $f$ has exactly two fixed points in $\mathbb{C}^3$. Suppose $f$ is reversible and has a fixed points $x_\pm$, then it is easy to see that $h(x_\pm)$ are also fixed points. In addition, $Df(x_\pm)$ is conjugate to $Df^{-1}(h(x_\pm))$.

Since there are two fixed points, either $h(x_\pm) = x_\pm$ or $h(x_\pm) = x_\mp$. Now the eigenvalues are invariant under a diffeomorphism, so in the first case the eigenvalues of $Df(x_+)$ must be the same as those of $Df^{-1}(x_+)$, and similarly for $x_-$. This can only happen when $t_\pm = s_\pm$ but this implies $x_+ = x_-$, which can not happen by assumption.

We conclude that $h(x_+) = x_-$. This implies that $t_\pm = s_\mp$, which gives $a = c$, or $Q(x,y) = Q(y,x)$.

The other direction is proved by a simple computation, as described above. ∎

Reversibility simplifies finding orbits of a map. Orbits that are invariant under $h$ are called symmetric, and as is easy to see, they must have points on the fixed set $Fix(h) = \{x \in \mathbb{R}^3 : h(x) = x\}$. In our case this is the line $x + z = -\eta, y = -\eta/2$, so a numerical algorithm for finding for symmetric orbits involves a one dimensional search.

Similarly, if the stable manifold of one of the fixed points intersects $Fix(h)$, then the intersection point is on a heteroclinic orbit, for suppose $x \in Fix(h) \cap W^s(x_+)$, then $x \in W^u(x_+)$, because $h(f^n(x)) = f^{-n}(h(x)) = f^{-n}(x)$, so

$$\lim_{n \to \infty} f^n(x) = x_+ \quad \Rightarrow \quad \lim_{n \to \infty} h(f^n(x)) = \lim_{n \to \infty} f^{-n}(x) = h(x_+) = x_-$$

Furthermore, suppose the stable manifold is two dimensional, and has the normal vector $\hat{n}$ at a point on $Fix(h)$, then $Dh(x)\hat{n}$ is the normal to the unstable manifold at this point. This implies that the curve of heteroclinic orbits is tangent to the direction $\hat{n} \times Dh(x)\hat{n}$.

## 5.4 Bounded orbits

For the Hénon map, it is well known that the set of bounded orbits is contained in a square. For the volume preserving case, we will show that an analogue of this result also holds providing the quadratic form $Q$ is positive definite:

**Theorem 5.1** *If $Q$ is positive definite then there is a $\kappa > 0$ such that all bounded orbits are contained in the cube $\{(x, y, z) : |x| \leq \kappa, |y| \leq \kappa, |z| \leq \kappa\}$. Moreover, points outside the cube go to infinity along the $+x$ axis as $t \to +\infty$ or the $-z$ axis as $t \to -\infty$*

**Proof:**

We start by writing the map in third difference form as

$$x_{t+1} = \alpha + \tau x_t - \sigma x_{t-1} + x_{t-2} + Q(x_t, x_{t-1}) .$$

Recall that a quadratic form $Q(x, y) = ax^2 + bxy + cy^2$ is positive definite *iff* $a > 0, c > 0$, and $d \equiv ac - b^2/4 > 0$. We will use the bounds obtained from completing the square:

$$\begin{aligned} Q(x, y) &= \frac{d}{c}x^2 + c(y + bx/2c)^2 \geq \frac{d}{c}x^2 , \\ &= (x + by/2a)^2 + \frac{d}{a}y^2 \geq \frac{d}{a}y^2 . \end{aligned}$$

There are three cases to consider, depending upon the relative sizes of $x_t$, $x_{t-1}$, and $x_{t-2}$:

- $|x_t| \geq max(|x_{t-1}|, |x_{t-2}|)$. The difference equation then gives

$$\begin{aligned} x_{t+1} &\geq Q(x_t, x_{t-1}) - |\alpha| - |\tau x_t| - |\sigma x_{t-1}| - |x_{t-2}| \\ &\geq \frac{d}{c}x_t^2 - (|\tau| + |\sigma| + 1)|x_t| - |\alpha| , \end{aligned}$$

Now since $d/c > 0$ there is a constant $\kappa_1 > 0$, depending on $\alpha, \tau, \sigma, a, b$, and $c$ such that when $|x_t| > \kappa_1$, we have

$$\frac{d}{c}x_t^2 - (|\tau| + |\sigma| + 1)|x_t| - |\alpha| > |x_t| ,$$

In this case, we have $x_{t+1} > |x_t|$. Noting that we then have $x_{t+1} > |x_t| \geq |x_{t-1}|$, we can recursively apply this result to show that the sequence

$$x_{t+k} > x_{t+k-1} > \ldots > |x_t| > \kappa_1$$

is monotone increasing. In fact, this sequence is unbounded; otherwise it would have a limit $x_t \to x^* > \kappa_1$, and this point would have to be a fixed point $x = y = z = x^*$ of the map. However, there are at most two such points, $x_\pm$, and a simple calculation shows that $\kappa_1 > x_\pm$, so both fixed point points are excluded.

- $|x_{t-2}| \geq max(|x_t|, |x_{t-1}|)$. Inverting the difference equation and shifting $t$ by one yields
$$x_{t-3} = x_t - \alpha - \tau x_{t-1} + \sigma x_{t-2} - Q(x_{t-1}, x_{t-2}) \ .$$

Thus we have
$$\begin{aligned} x_{t-3} &\leq -\frac{d}{a} x_{t-2}^2 + |x_t| + |\alpha| + |\tau x_{t-1}| + |\sigma x_{t-2}| \\ &\leq -\frac{d}{a} x_{t-2}^2 + (|\tau| + |\sigma| + 1)|x_{t-2}| + |\alpha| \\ &< -|x_{t-2}| \ , \end{aligned}$$

when $|x_{t-2}| > \kappa_2$, for a constant $\kappa_2$ chosen as before, but with $d/c$ replacing $d/a$. This implies that the sequence $x_{t-k} < x_{t-k+1} < \ldots < -|x_{t-2}|$ is monotone decreasing, negative, and unbounded.

- $|x_{t-1}| \geq max(|x_t|, |x_{t-2}|)$. In this case we will see that the orbit is unbounded in both directions of time. For the forward direction, note that
$$\begin{aligned} x_{t+1} &\geq \frac{d}{a} x_{t-1}^2 - (|\tau| + |\sigma| + 1)|x_{t-1}| - |\alpha| \\ &> |x_{t-1}| \ , \end{aligned}$$

when $|x_{t-1}| > \kappa_2$. Thus $x_{t+1} > |x_{t-1}| \geq |x_t|$, which is the situation covered by (i), and we get a monotone increasing sequence, providing $x_{t+1} > \kappa_1$. Alternatively, note that
$$\begin{aligned} x_{t-3} &\leq -\frac{d}{c} x_{t-1}^2 + (|\tau| + |\sigma| + 1)|x_{t-1}| + |\alpha| \\ &< -|x_{t-1}| \ , \end{aligned}$$

when $|x_{t-1}| > \kappa_1$. This gives $x_{t-3} < -|x_{t-1}| < -|x_{t-2}|$, so we are in the situation covered by (ii), which implies that the sequence approaches $-\infty$ providing $|x_{t-3}| > \kappa_2$

In conclusion, we have shown that an orbit is unbounded either as $t \to \pm\infty$ providing it contains a point $x_t$ such that $|x_t| > max(\kappa_1, \kappa_2) \equiv \kappa$. Note that $\kappa_1$ is a monotone decreasing function of $d/c$, so therefore we can define $\kappa$ by using $d/max(a, c)$:

$$\kappa = \frac{max(a,c)}{2d}\left(|\tau| + |\sigma| + 2 + \sqrt{(|\tau| + |\sigma| + 2)^2 + 4\frac{|\alpha|}{d} max(a,c)}\right)$$

Finally, we investigate the asymptotic direction of an unbounded orbit. Recall that $y_t = x_{t-1}$ and $z_t = x_{t-2}$. Suppose that $|x_t| > |y_t| > |z_t| > \kappa$, then each of the variables is eventually positive, so the orbit moves to infinity in the positive octant. Moreover, once all components are positive, we have

$$\begin{aligned}\frac{x_{t+1}}{x_t} &= \frac{Q(x_t,y_t)}{x_t} + \frac{\alpha + \tau x_t - \sigma y_t + z_t}{x_t} \\ &\geq \frac{d}{c}x_t - (|\tau|+|\sigma|+1) - \frac{|\alpha|}{x_t} \to \infty \ .\end{aligned}$$

So the ratios $y_t/x_t = x_{t-1}/x_t$ and $z_t/y_t = x_{t-2}/x_{t-1}$ go to zero, and the orbit approaches the positive $x$-axis as $t \to \infty$. Similarly if $|z_t| > |y_t| > |x_t| > \kappa$, then eventually all the components are negative, and so the orbit moves to infinity in the negative octant as $t \to -\infty$. Once all components are negative, we have

$$\begin{aligned}\frac{z_{t-1}}{z_t} &= \frac{x_{t-3}}{x_{t-2}} \\ &= -\frac{Q(y_t,z_t)}{z_t} + \frac{x_t - \alpha - \tau y_t + \sigma z_t}{z_t} \\ &\leq -\frac{d}{a}|z_t| + (|\tau|+|\sigma|+1) + \frac{|\alpha|}{|z_t|} \to -\infty \ .\end{aligned}$$

This implies the orbit moves to $\infty$ along the negative $z$ axis. ∎

# 6   Conclusions

We have studied a family of volume preserving maps with the property that all entries are quadratic polynomials. We showed that these conditions imply that such maps are polynomial diffeomorphisms. Then we restricted ourselves to quadratic maps whose inverse is also quadratic. The class of maps studied is related to an old conjecture about polynomial maps called the Jacobian conjecture.

A definition of quadratic shears was introduced and a characterization was given in general. A further characterization in three space was applied to find a normal form for the family. In three-space, the form of the generic case is similar in form to the area preserving Hénon map and, generically, the map has two fixed points that can be either *type A* or *type B*.

In addition, using our definition of quadratic shear and its characterization, we were able to give a simpler proof of a theorem of Moser classifying quadratic symplectic maps.

Once we have found the normal form for the quadratic volume preserving map, there remain many enticing open problems. For example, we plan further computations to visualize the stable and unstable manifolds of the fixed points. Often these manifolds intersect, enclosing a ball; however, this is not guaranteed. Moreover, the heteroclinic intersections, which are generically curves, can fall in many homotopically distinct classes. We suspect that there are bifurcations between these classes, and that which occurs will depend, for example, on the complex phase of the eigenvalue of the associated fixed point. Heteroclinic orbits can be found most easily for the reversible case, as an intersections should occur on the fixed set of the reversor

Another problem of interest is to obtain a characterization of quadratic shears in higher dimensions similar to the one we obtained in three dimensions. At this point, normal forms could be obtained using techniques similar to the current paper.

Finally, as we discussed in the introduction, one of our main motivations for characterizing the quadratic volume preserving maps is to study transport. If the two fixed points have disparate types, and their two dimensional manifolds intersect on a circle, then transport can be localized to "lobes" similar to the two dimensional case [18]. However, as Fig. 2 shows, the intersections can be curves that spiral from one fixed point to the other. We plan to characterize transport for such cases. The existence, for the definite case, of a cube containing the bounded orbits (cf. theorem 5.1) will prove useful in this study.

# 7 Figures

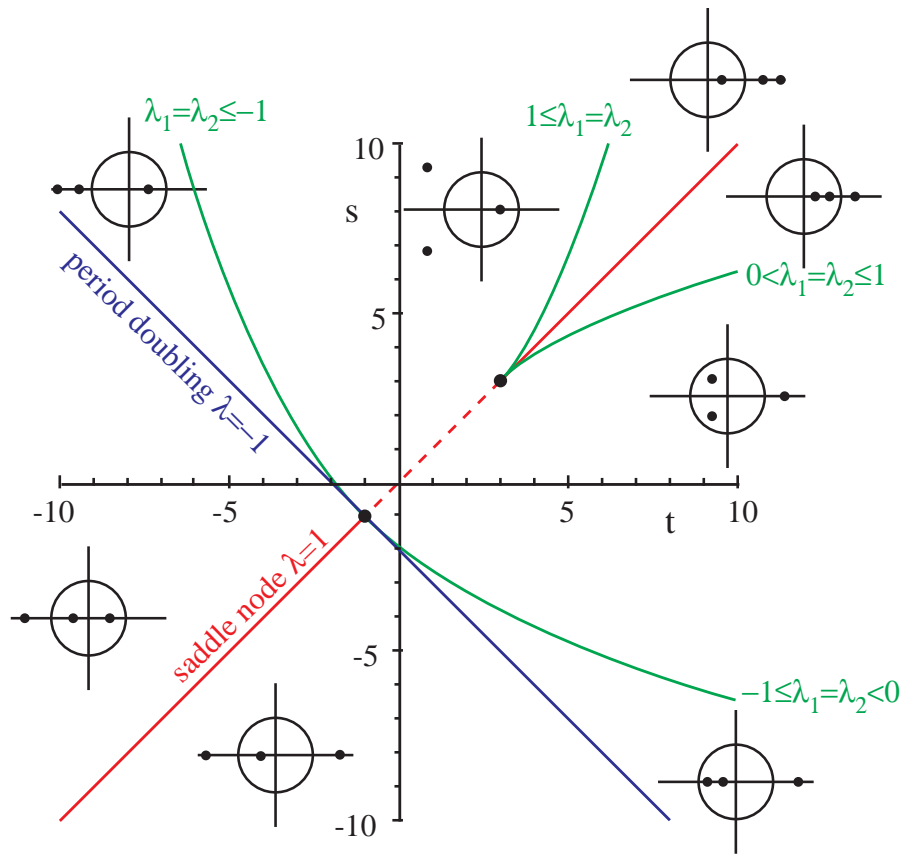

Figure 1: General stability diagram for a volume preserving map.

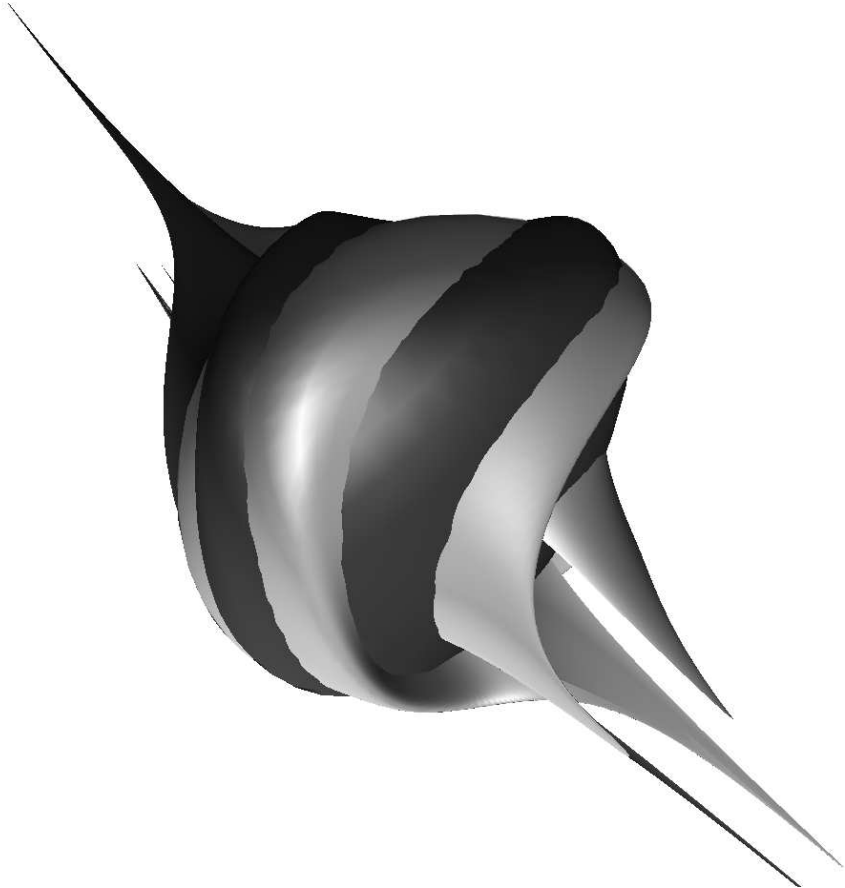

Figure 2: Two dimensional stable and unstable manifolds for the parameters $a = c = 0.5, b = 0.0, \alpha = 0.0, \tau = -0.3$ and $\sigma = 0.0$.

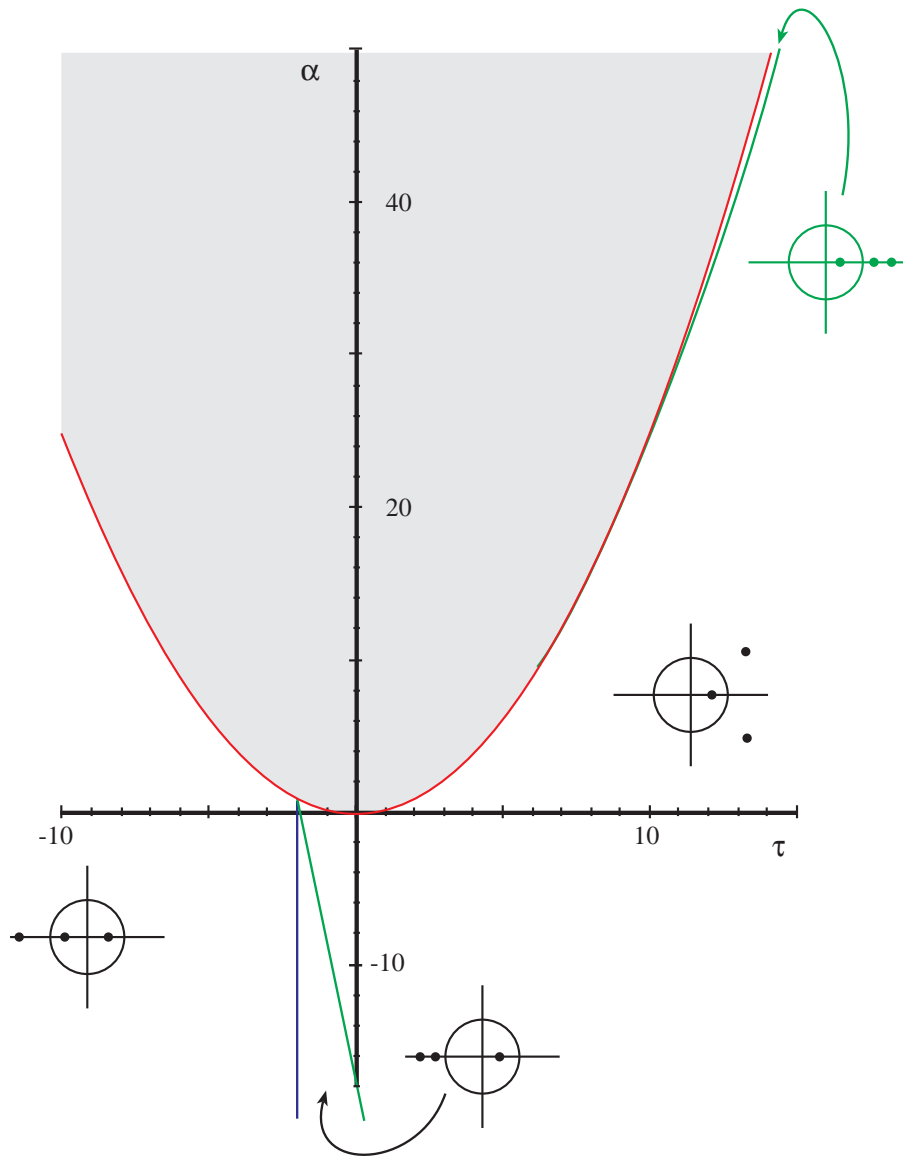

Figure 3: Stability Diagram for the reversible case $a + b + c = 1$, $a = c$ and $\sigma = 0$. There are no fixed points in the shaded region. The complex eigenvalues for the fixed point $x_+$ are shown. Those of $x_-$ are reciprocal to those of $x_+$.

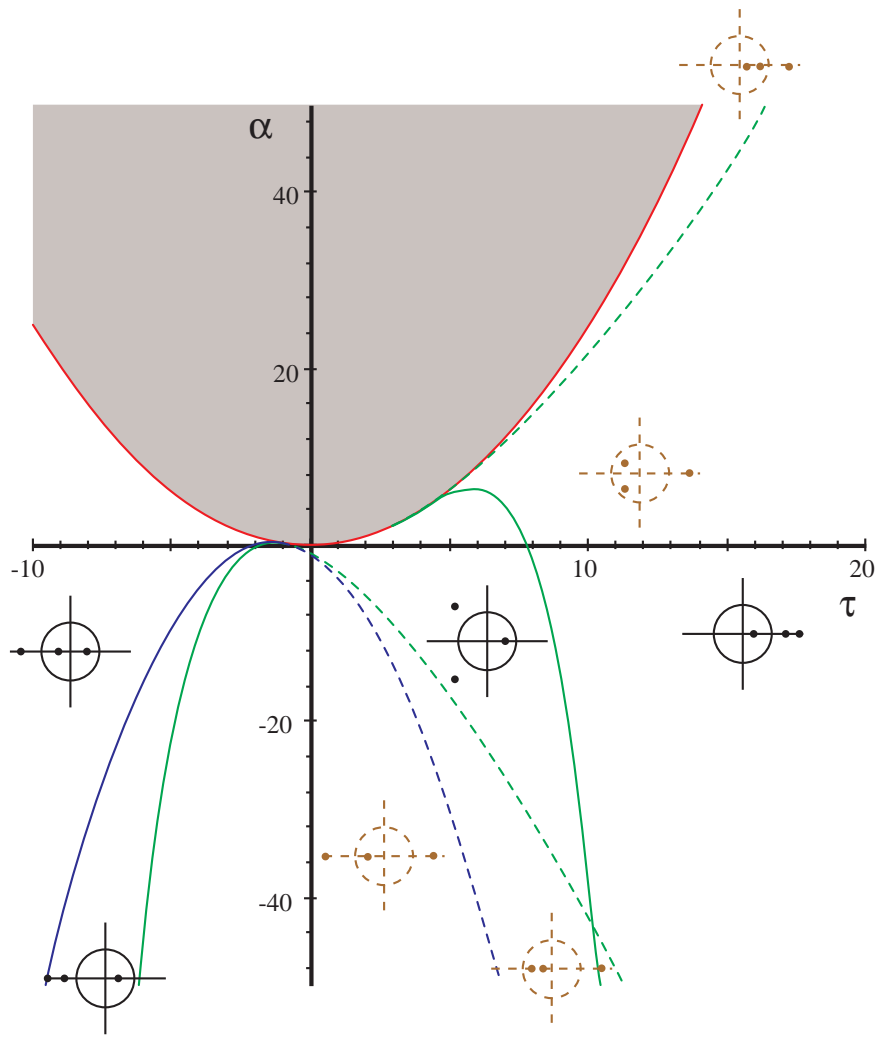

Figure 4: Stability Diagram for $a = -1/2, b = 1, c = 1/2, \sigma = 0$. Solid lines correspond to changes in the stability of $x_+$, while dashed lines to $x_-$.